\documentclass[3p,sort&compress]{elsarticle}

\pdfoutput=1

\usepackage[utf8]{inputenc}         
\usepackage[T1]{fontenc}            
\usepackage[british]{babel}         
\usepackage{amsmath}                
\usepackage{amsfonts}               
\usepackage{hyperref}               
\usepackage{graphicx}               
\usepackage{microtype}              

\newcommand{\mytitle}{Control-based continuation: bifurcation and stability
  analysis for physical experiments}

\hypersetup{colorlinks=true, allcolors=blue}
\makeatletter
\AtBeginDocument{
  \hypersetup{pdftitle={\mytitle}, pdfauthor={David A.W. Barton}} 
}
\makeatother

\graphicspath{{.}{./figures/}}

\newcommand{\eref}[1]{Eq.~(\ref{#1})}

\newcommand{\fref}[1]{Fig.~\ref{#1}}
\newcommand{\Fref}[1]{Figure~\ref{#1}}
\newcommand{\sref}[1]{Sec.~\ref{#1}}

\DeclareMathOperator{\md}{d}

\makeatletter
\def\diff{\@ifnextchar[{\@diffwith}{\@diffwithout}}
\def\@diffwith[#1]#2#3{\frac{\md^{#1} #2}{\md #3^{#1}}}
\def\@diffwithout#1#2{\frac{\md #1}{\md #2}}
\makeatother

\def\real{\mathbb{R}}

\def\mm{\,\text{mm}}

\def\ms{\,\text{m}\text{s}^{-1}}
\def\ms2{\,\text{m}\text{s}^{-2}}
\def\Hz{\,\text{Hz}}
\def\kHz{\,\text{kHz}}

\def\vec#1{\ensuremath{\mathbf{#1}}}

\usepackage[svgnames]{xcolor}

\def\newtext#1{{#1}}

\begin{document}


\begin{frontmatter}


  \title{\mytitle}
  \author[ENM]{David~A.W.~Barton\corref{DB}}
  \ead{david.barton@bristol.ac.uk}
  \address[ENM]{Department of Engineering Mathematics, University of Bristol}
  \cortext[DB]{Corresponding author}
  

  \begin{abstract}
    Control-based continuation is technique for tracking the solutions and
    bifurcations of nonlinear experiments. The basic idea is to apply the method
    of numerical continuation to a feedback-controlled physical
    experiment. Since in an experiment it is not (generally) possible to set the
    state of the system directly, the control target is used as a proxy for the
    state. The challenge then becomes to determine a control target such that
    the control is non-invasive, that is, it stabilises a steady-state (or
    periodic orbit) of the original open-loop experiment without altering it
    otherwise.

    Once implemented, control-based continuation enables the systematic
    investigation of the bifurcation structure of a physical system, much like
    if it was numerical model. However, stability information (and hence
    bifurcation detection and classification) is not readily available due to
    the presence of stabilising feedback control. This paper uses methods from
    the system identification community to extract stability information in the
    form of Floquet multipliers from the closed-loop experiment, thus enabling
    the direct detection of bifurcations. In particular, it is shown that a
    periodic auto-regressive model with exogenous inputs (ARX)
    can be constructed that approximates the time-varying linearisation of the
    experiment around a particular periodic orbit. This method is demonstrated
    using a physical nonlinear tuned mass damper.
  \end{abstract}


\end{frontmatter}


\section{Introduction}\label{sec:intro}


Control-based continuation is a systematic method for performing bifurcation
studies on physical experiments. Based on modern feedback control schemes it
enables dynamical phenomena to be detected and tracked as system parameters are
varied in a similar manner to how nonlinear numerical models can be investigated
using numerical continuation. Control-based continuation was originally
developed as an extension to Pyragas' time-delayed feedback
control~\cite{Pyragas1992, Pyragas2001, Pyragas2002} to make it more robust and
suitable for parameter studies~\cite{SieberKrauskopf2008}, though its current
implementation contains no elements of time-delayed feedback.


The use of feedback control for the investigation of nonlinear systems is not
new; in addition to time-delayed feedback, methods such as OGY
control~\cite{OttGrebogiYorke1990} have previously been employed to provide
non-invasive control to stabilise unstable orbits and investigate dynamical
phenomena. Indeed, previous authors have gone as far as to implement such
control schemes within parameter continuation studies to numerically simulate
particular experiments such as atomic force
microscopes~\cite{MisraDankowiczPaul2008} or reaction
kinetics~\cite{SiettosKevrekidisMaroudas2004}. Control-based continuation goes
beyond these particular methods to allow the use of almost any feedback control
scheme and, as such, it is a general purpose tool applicable to a wide range of
physical experiments.


Control-based continuation has been successfully applied to a range of
experiments including a parametrically excited
pendulum~\cite{SieberGonzalez-BuelgaNeildEtAl2008}, nonlinear energy
harvesters~\cite{BartonBurrow2011, BartonMannBurrow2012} and a bilinear
oscillator~\cite{BureauSchilderEtAl2013, BureauSchilderEtAl2014}. In each case,
periodic orbits have been tracked through instabilities such as saddle-node
bifurcations (folds) thus revealing a great deal of dynamical information about
the system in question, including the location of a codimension-2 cusp
bifurcation in one case~\cite{BartonSieber2013}. As well as bifurcations, other
dynamic features such as backbone curves can also be tracked with control-based
continuation~\cite{Renson2015}.



Although the basic scheme for control-based continuation is well established (an
overview is provided in \sref{sec:cbc}), it lacks many of the features of
standard numerical continuation schemes such as bifurcation detection. Only
saddle-node bifurcations (folds) can be detected readily and that is because
they are geometric features in the solution surface. Bifurcations such as
period-doubling bifurcations are not geometric features and can go undetected
due to the stabilising affect of the feedback controller. \newtext{Similarly,
  the inclusion of the feedback controller means that methods for calculating
  eigenvalues/Floquet multipliers and basins of attraction from experiments such
  as~\cite{MurphyBaylyEtAl1994, VirginToddEtAl1998, KozinskyPostmaEtAl2007} are
  not helpful since they indicate the stability of the closed-loop system rather
  than the open-loop system.}

In this paper we consider only periodically forced systems and hence study
periodic orbits, though there is no reason that the methods developed should not
be applicable to autonomous systems as well. In \sref{sec:linear}, we propose a
method for calculating the stability (the Floquet multipliers and associated
stable and unstable eigendirections) based on the estimation of a local
linearisation around a stabilised periodic orbit. We demonstrate the
effectiveness of this method in \sref{sec:apparatus} and \sref{sec:results} by
applying it to a (physical) nonlinear mass-spring-damper system where the
nonlinearity is geometric in nature --- the springs are mounted perpendicular to
the direction of motion.


\section{Control-based continuation}\label{sec:cbc}


Numerical continuation is a path following method used to track solution
branches as parameters of the system in question are varied. In a nonlinear
system, these solution branches can encounter bifurcations at particular
parameter values which results in a qualitative change in the dynamics of the
system. Numerical continuation enables these bifurcations to be detected and
tracked in turn.  It is typically applied to differential equation models but it
can be used more widely, for example on finite element models.


At a basic level, numerical continuation tracks the solutions of an arbitrary
nonlinear function, a \emph{zero problem} given by
\begin{equation}
  \label{eqn:f}
  f(x, \lambda) = 0, \qquad f:\real^n\times\real^p\rightarrow\real^m
\end{equation}
where $x$ is the system state and $\lambda$ is the system
parameter(s). \newtext{A common example of this is tracking the equilibria of an
  ordinary differential equation with respect to a single parameter. In this
  case $n=m$ and $p=1$ in \eref{eqn:f}, that is, $f=0$ defines a one-dimensional
  curve. Alternatively the function $f$ can arise from the discretisation of
  a periodic orbit.} Numerical continuation works in a predictor-corrector
fashion; at each step a new solution $\tilde{x}$ is predicted from previously
determined solutions and then the solution is corrected using a nonlinear root
finder applied to the function $f$ (typically a Newton iteration). The use of a
nonlinear root finder means that the stability or instability of solutions is
unimportant. In certain circumstances (for example, near a fold or saddle-node
bifurcation) the function $f$ must be augmented with an additional equation ---
the pseudo-arclength equation --- which enables the numerical continuation
scheme to track solution curves that double back on themselves. In these
circumstances, without the pseudo-arclength equation the correction step for a
fixed set of parameter values $\lambda$ will fail since no solution exists. For
extensive information and guidance on numerical continuation see the
textbooks~\cite{Seydel2010, Kuznetsov1998}. Numerical software is also readily
available in the form of \textsc{CoCo}~\cite{DankowiczSchilder2013} and
AUTO~\cite{DoedelOldeman2012} amongst others.

Control-based continuation is a means for defining a zero-problem based on the
outputs of a physical experiment, thus enabling numerical continuation to be
applied directly without the need for a mathematical model. To do this there are
two key challenges to overcome. 1) In general, it is not possible to set the
state $x$ of the physical system and so it is not possible to evaluate $f$ at
arbitrary points. 2) The physical system must remain around a stable operating
point while the experiment is running. While a numerical model going unstable
might prove to be a mild annoyance, a physical system going unstable can prove
dangerous.

In order to overcome these challenges, a feedback controller is used to
stabilise the system and the control target (or reference signal) acts as a
proxy for the system state. \newtext{The feedback controller takes the form
\begin{equation}
  u(t) = g(x^*(t) - x(t))
  \label{eqn:control}
\end{equation}
where $x^*(t)$ is the control target and $g$ is a suitable control law such as
proportional-derivative (PD) control (as used in this paper) where
\begin{equation}
  u(t) = K_p(x^*(t) - x(t)) + K_d(\dot x^*(t) - \dot x(t)).
\end{equation}
For the method outlined in this paper, the choice of control law is at the
discretion of the user; any suitable stabilising feedback control scheme can be
used. The challenge here is to devise a scheme for embedding the feedback
control within the numerical continuation such that the controller becomes}
non-invasive, that is, the controller does not affect the locations of any
invariant sets in the experiment such as equilibria or period orbits. This
requirement for non-invasiveness defines the zero problem; a control target must
be chosen such that the control action
\begin{equation}
  u(t)\equiv0.
  \label{eqn:control_zero}
\end{equation}
  
In this paper, we consider the case of a periodically forced experiment
\newtext{with forcing frequency $\omega$} and, as such, only consider periodic
motions. In this case it is appropriate to consider a Fourier discretisation of
\eref{eqn:control_zero} and so find the coefficients of the Fourier series of
the control target \newtext{$x^*(t)=A^{x^*}_0/2 +
  \sum_{j=1}^{m}A^{x^*}_j\cos(j\omega t) + B^{x^*}_j\sin(j\omega t)$ such that
  \eref{eqn:control_zero} is satisfied}. (In other circumstances different
discretisations may be appropriate.) \newtext{In this case the control action
  $u$ has a Fourier series representation given by}
\begin{equation}
  u(t) = \frac{A^u_0}{2} + \sum_{j=1}^{m}A^u_j\cos(j\omega t) + B^u_j\sin(j\omega t).
  \label{eqn:u_Fourier}
\end{equation}
\newtext{where the Fourier coefficients $A^u_j$ and $B^u_j$ are derived directly
  from the measured control action~\eref{eqn:control}, that is,
\begin{subequations}
  \label{eqn:u_fourier}
  \begin{gather}
    A^u_j =
    \frac{\omega}{\pi}\int_0^{\frac{2\pi}{\omega}}g(x^*(t)-x(t))\cos(j\omega
    t)\md t,\quad\text{for}\ j=0, 1, 2, \ldots\\
    B^u_j =
    \frac{\omega}{\pi}\int_0^{\frac{2\pi}{\omega}}g(x^*(t)-x(t))\sin(j\omega
    t)\md t,\quad\text{for}\ j=1, 2, \ldots.
  \end{gather}
\end{subequations}
Hence the discretised zero problem is defined as
\begin{equation}
  0 = A^u_j, \quad 0=B^u_j \quad \forall j.
  \label{eqn:cbc_zero}
\end{equation}}

To solve \eref{eqn:cbc_zero} standard root finding algorithms can be used,
however, any required derivatives must be estimated using finite differences
from experimental data \newtext{after adjusting the experiment
  inputs}. Consequently, gradient-based methods can be slow despite their good
convergence rates. In previous publications a Newton-Broyden method, which
avoids recomputing derivative information for successive iterates, has proven
effective~\cite{SieberGonzalez-BuelgaNeildEtAl2008, BartonMannBurrow2012,
  BureauSchilderEtAl2013}.

In this paper, since the control acts through the same mechanism as the forcing,
we are able to use a quicker method which exploits the fact that we are
performing a parameter study in the forcing
amplitude~\cite{BartonSieber2013}. Consider the case where the total input to
the system is given by
\begin{equation}
  i(t) = p(t) + u(t)
  \label{eqn:total_input}
\end{equation}
where $p(t)$ is the forcing signal \newtext{and $u(t)$ is the control
  action. Furthermore, we consider} the case of sinusoidal forcing where
\begin{equation}
  p(t) = a\cos(\omega t) + b\sin(\omega t).
\end{equation}

For an arbitrary control target \newtext{$u^*(t)$ and forcing input $p(t)$}, the
Fourier coefficients of \eref{eqn:u_Fourier} will be non-zero. However, the
contribution in the fundamental mode (coefficients $A^u_1$ and $B^u_1$) can be
lumped together with the coefficients $a$ and $b$ of the forcing term giving a
new effective forcing amplitude of
\begin{equation}
  \Gamma = \sqrt{\left(a+A^u_1\right)^2 + \left(b+B^u_1\right)^2}.
\end{equation}
\newtext{Hence, once the higher Fourier modes of the control action $u(t)$ are
  set to zero (as described below), the total input to the system will be $i(t)
  = \Gamma\cos(\omega t + \phi)$ (the phase $\phi$ is unimportant since the
  system is time invariant).}  In essence, instead of setting the forcing
amplitude and trying to calculate the correct corresponding control target, we
set the control target and measure the corresponding forcing amplitude.

Though this procedure leaves the $A^u_0$ and the higher harmonics untouched, the
corresponding control target coefficients required to set the control action to
zero can be quickly determined using a fixed-point iteration. In an iterative
manner the remaining coefficients of the control target $x^*(t)$ are simply set
equal to the measured coefficients of the system response $x(t)$, the iteration
finishes when the response and the control target remain equal for a certain
period of time. For the system described below, this takes a single iteration.

A fixed-point iteration cannot be applied to the coefficients of the fundamental
mode since, generically, instabilities in the system will manifest in the
fundamental mode.

\newtext{It is important to emphasise that this procedure does not depend on the
  specifics of the control law used in the experiment. All that is required is
  that smooth changes in the control target $x^*(t)$ result in smooth changes in
  the Fourier coefficients~\eref{eqn:u_fourier}. As such, this method is
  convenient in realistic settings where the control law is more sophisticated
  than simple PD control and filtering of the signal is required; these effects
  are simply captured in the control action $u(t)$ and fed
  into~\eref{eqn:u_fourier}.}


\section{Identification of a linearisation}\label{sec:linear}

For a typical ordinary differential equation model, \newtext{the right-hand side
  of which is given by $h(x(t))$,} the Floquet multipliers and hence stability
of a periodic orbit are determined by integrating over one period the first
variational equation
\begin{equation}
  \diff y t = A(\hat x(t))y(t)  \label{eqn:variational}
\end{equation}
where $y(t)$ represents the deviation from \newtext{a predetermined periodic
  orbit $\hat x(t)$ and the Jacobian matrix $A$ is calculated from the
  derivatives of $h$ with respect to $x$ evaluated along the periodic orbit
  given by $\hat x(t)$ (as such $A$ is a time varying quantity).}

In the context of control-based continuation, even determining whether an orbit
is stable or unstable is problematic due to the presence of stabilising feedback
control. In \cite{BureauSchilderEtAl2014} a number of measures are suggested to
overcome this problem but all require turning off the feedback control for a
period of time; in many situations this is not desirable as damage could be
caused to the experiment or even the experimenter. As such, here we attempt to
fit a time-varying linearisation using techniques from the system identification
community. From the fitted time-varying linearisation we are able to calculate
the corresponding Floquet multipliers and hence the stability properties of the
periodic orbit.

\newtext{We start by assuming that the experiment of interest is undergoing a
  periodic motion $\hat x(t)$ for a given forcing input $\hat i(t) = \hat p(t) +
  \hat u(t)$ (cf.~\eref{eqn:total_input}). In order to generate data with which
  to fit a time-varying linearisation, the system is perturbed using filtered
  Gaussian white noise $\eta(t)$ such that the total input to the experiment is
  \begin{equation}
    i(t) = \hat i(t) + \eta(t) + u(t) - \hat u(t).
  \end{equation}
  The additional $u(t) - \hat u(t)$ term arises due to the presence of the
  feedback controller acting against the applied perturbation. (Details of
  $\eta(t)$ are below.) Finally, we define the perturbed system response
  \begin{equation}
    y(t) = x(t) - \hat x(t)
  \end{equation}
  and, similarly, the perturbed system input
  \begin{equation}
    k(t) = i(t) - \hat i(t).
  \end{equation}
}
  
\begin{figure}
  \centering
  \includegraphics{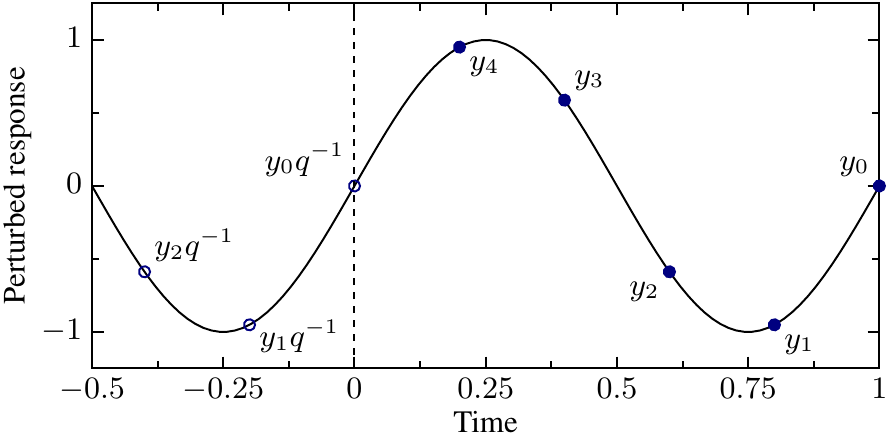}
  \caption{An example discretisation of the perturbed response of a periodic
    orbit with $m=5$ and $n=2$; the period of the orbit is normalised to
    $1$. Here $\vec y = [y_0,y_1,y_2,y_3,y_4]$. Since $n=2$, $y_0$ is purely a
    function of $[y_1, y_2]$ and $[k_0,k_1,k_2]$. Similarly, $y_4$ is purely a
    function of $[y_0q^{-1}, y_1q^{-1}]$ and $[k_4,k_0q^{-1},k_1q^{-1}]$, where
    $q^{-1}$ is the backward shift operator. Thus the linearised model
    \eref{eqn:arx} allows the construction of a linear period map from $\vec
    yq^{-1}$ to $\vec y$.}
  \label{fig:sysid}
\end{figure}

Rather than try to fit a continuous time model of the response to perturbations
such as \eref{eqn:variational}, which requires the estimation of derivatives
from experimental data, we instead fit the coefficients of a discrete-time
multiple-input multiple-output (MIMO) auto-regressive model with exogenous
inputs (ARX) of the form
\begin{equation}
  B(q^{-1})\vec y(T) = A(q^{-1})\vec k(T) + \vec e(T)
  \label{eqn:arx}
\end{equation}
where $\vec y(T) = [y(T-i/m)]_{i=0\ldots m-1}$, $\vec k(T) =
[k(T-i/m)]_{i=0\ldots m-1}$ and $\vec e(T) = [e(T-i/m)]_{i=0\ldots m-1}$ are
vectors of data points of the perturbed system response, the perturbed system
input (due to the control action and an additional random perturbation) and the
model error respectively, sampled across a single period of oscillation and
$q^{-1}$ is the backward shift
operator~\cite[Sec.~6.2]{SoderstromStoica1989}. Thus $\vec y(T)$ corresponds to
a discretisation of the perturbed system response using $m$ points over the
period; this discretisation is illustrated in \fref{fig:sysid}. $B(q^{-1})$ and
$A(q^{-1})$ are square $m\times m$ matrices of polynomials in $q^{-1}$; here we
restrict the polynomials to being first order (at most). Thus \eref{eqn:arx}
acts as a period map with all the dynamics of interest encoded in
$B(q^{-1})$. (The matrix $A(q^{-1})$ though required for system identification
are not used to infer stability.)

\newtext{The model error $\vec e(T)$ is not measured directly but is minimised
  by the particular system identification method used on~\eref{eqn:arx}.}

\newtext{ARX models are used extensively in the system identification and time
  series analysis communities. Their simplicity has seen them applied to a wide
  range of topics. They are particularly appropriate in situations where
  discretely sampled data is available as they avoid the need of estimating
  derivatives. For more information see the textbook by
  Hamilton~\cite{Hamilton1994} or one of the many other books on this topic.}

At each point we assume linear observability to state that the next point in the
time evolution of the linearisation is determined entirely by the previous $n$
points alone. We assume that $n < m$ (increase the value of $m$ as appropriate)
to obtain a banded matrix structure for $B(q^{-1})$ of
\begin{equation}
  B(q^{-1}) = 
  \begin{bmatrix}
    1 & b_{1,1} & \cdots & b_{1,n-1} & b_{1,n} & 0 & 0 & \cdots & 0\\
    0 & 1 & \cdots & b_{2,n-2} & b_{2,n-1} & b_{2,n} & 0 & \cdots & 0\\
    \ddots & \ddots & \ddots & \ddots & \ddots & \ddots & \ddots & \ddots &
    \ddots\\
    b_{m,1}q^{-1} & b_{m,2}q^{-1} & \cdots & b_{m,n}q^{-1} & 0 & 0 & 0 & \cdots & 1
  \end{bmatrix},
\end{equation}
and for $A(q^{-1})$
\begin{equation}
  A(q^{-1}) = 
  \begin{bmatrix}
    a_{1,0} & a_{1,1} & \cdots & a_{1,n-1} & a_{1,n} & 0 & 0 & \cdots & 0\\
    0 & a_{2,0} & \cdots & a_{2,n-2} & a_{2,n-1} & a_{2,n} & 0 & \cdots & 0\\
    \ddots & \ddots & \ddots & \ddots & \ddots & \ddots & \ddots & \ddots &
    \ddots\\
    a_{m,1}q^{-1} & a_{m,2}q^{-1} & \cdots & a_{m,n}q^{-1} & 0 & 0 & 0 & \cdots & a_{m,0}
  \end{bmatrix}.
\end{equation}
Thus the system identification procedure must identify the $m(2n+1)$
coefficients \newtext{($a_{i,j}$ and $b_{i,j}$) within these matrices} to fully identify the
linear model. The optimal values for $m$ and $n$ can be estimated using the
Akaike information criterion (AIC)~\cite{BrockwellDavis1991} or similar. Note
that while increasing $n$ increases the data requirements for successful system
identification, increasing $m$ does not since more information is taken from the
existing time series.

Since the experiment is operating in closed loop, the number of available system
identification methods is somewhat limited. In this paper we make use of the
direct method for closed-loop identification due to its simplicity \newtext{and
  so identify the unknown parameters of \eref{eqn:arx} with least squares, thus
  minimising the sum-of-squares of the model error $\vec e(T)$.} However, other
methods such as joint input-output identification can be used if
required~\cite{SoderstromStoica1989, Ljung1999}. In order to provide
sufficiently informative results for system identification purposes, the random
perturbation \newtext{$\eta(t)$} should have a sufficiently broadband
spectrum. However, to minimise extraneous noise (that is, random perturbations
which are not captured by the discretisation) the bandwidth of the random
perturbation should be less than the Nyquist frequency corresponding to the
discretisation in \eref{eqn:arx}. \newtext{For the purposes of this paper,
  $\eta(t)$ is generated by passing Gaussian white noise through a 6th order
  Butterworth filter with a cut-off frequency of $10\,$Hz.}

\newtext{When there is a significant amount of measurement noise or unmeasurable
  random disturbances, a non-trivial noise model is required in \eref{eqn:arx},
  giving rise to a moving-average (MA) term. In this case the the unknown
  coefficients of the resulting ARMAX model} must be estimated using a method
such as the prediction error method (PEM) since the straightforward use of
linear least-squares will result in bias~\cite[Chapter
  10]{SoderstromStoica1989}. However, linear least-squares provides a quick and
effective way of starting the iterative PEM optimisation.

Once a linearised model has been identified the Floquet multipliers of the
periodic orbit can be determined from the matrix
$B(q^{-1})$. \newtext{Specifically, we seek to determine the monodromy matrix
  $M$ such that $\vec y(T) = Mq^{-1}\vec y(T)$ with $\vec k(T)\equiv 0$ and
  $\vec e(T)\equiv 0$, that is, we seek a linear mapping which takes one period
  of data points and returns the following period of data points subject to no
  disturbance to the system input. For the method described here $M$}
corresponds to the first $n$ rows and $n$ columns of the matrix given by
\begin{equation}
  B^{-1}(0)(B(1) - B(0)),
\end{equation}
and the Floquet multipliers of the periodic orbit are the eigenvalues of $M$. In
addition to the Floquet multipliers, the eigenvectors of $M$ correspond to the
stable and unstable directions of the periodic orbit at a particular point in
the oscillation.


\section{Experimental apparatus}\label{sec:apparatus}

\begin{figure}
  \centering
  \includegraphics{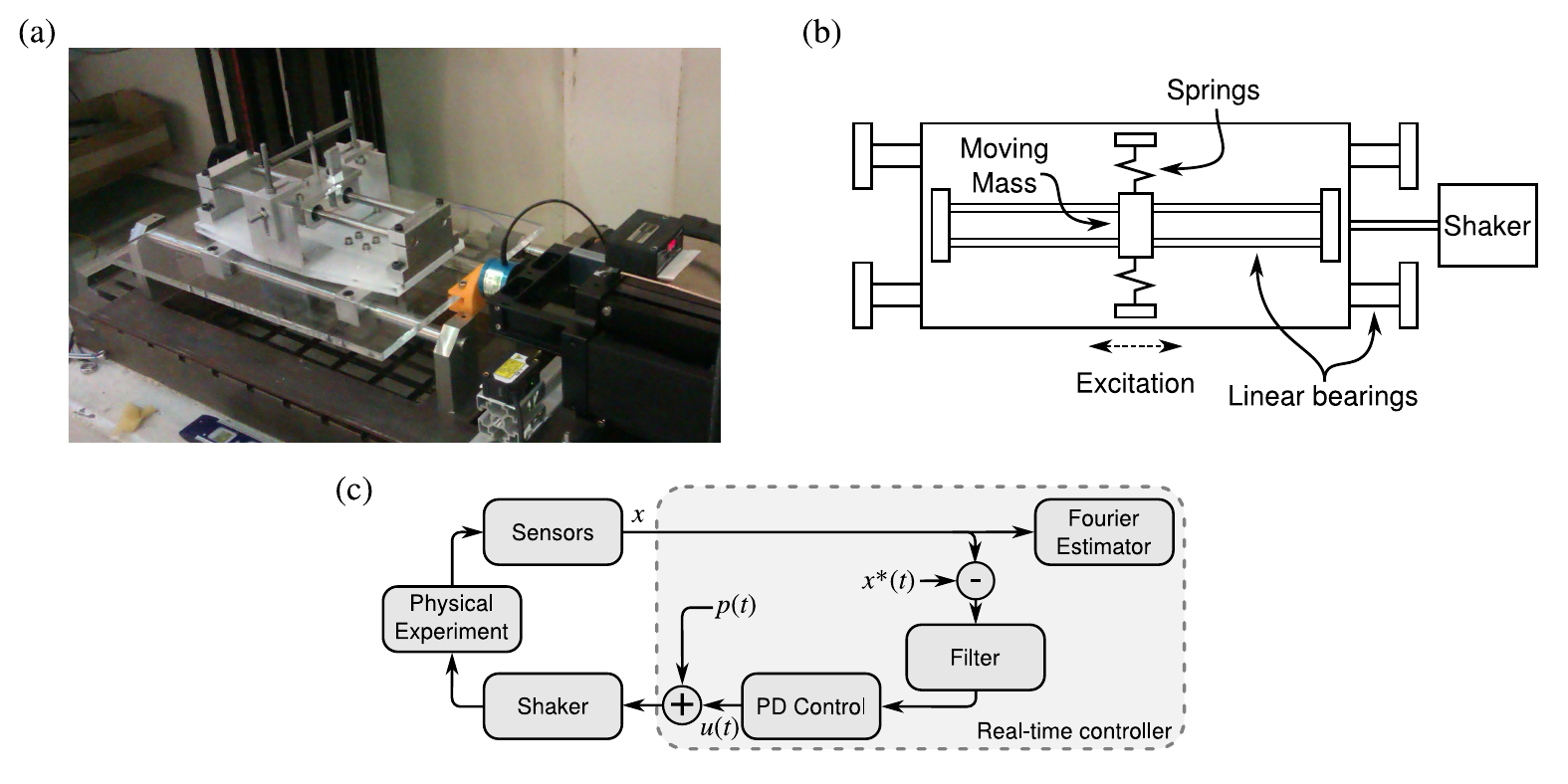}
  \caption{Panel (a) shows photograph of the nonlinear tuned mass damper (NTMD)
    used to test the effectiveness of the methodology outlined in this
    paper. Panel (b) shows a schematic of the NTMD with the springs mounted
    perpendicular to the direction of motion which results in a geometric
    nonlinearity. Panel (c) shows a schematic of the overall experimental rig;
    the feedback-control loop and a limited amount of signal processing are
    implemented in real-time while the numerical continuation routines are
    implemented off-line \newtext{(that is, there are no time constraints on the
      computations for the continuation; the experiment will simply continue
      running until new input parameters are available).}}
  \label{fig:exp}
\end{figure}

To test the effectiveness of this methodology outlined in this paper we apply
control-based continuation to a nonlinear tuned mass damper (NTMD) similar to
the one described in \cite{AlexanderSchilder2009}. The NTMD consists of a mass
able to move horizontally on a low friction bearing system while restrained by
two springs that are mounted perpendicularly to the direction of motion, thus
providing a geometric nonlinearity. The NTMD is then excited at the base. This
configuration results in a hardening spring-type characteristic. A photograph of
the experiment is shown in \fref{fig:exp}(a) along with a schematic of the
experiment in \fref{fig:exp}(b).

The details of the actuation and measurement equipment are as follows. The NTMD
is excited using an APS-113 long-stroke electrodynamic shaker in current control
mode using a Maxon ADS-50/10-4QDC motor controller. Typical base displacements
are sinusoidal with a frequency ranging from $2.2$--$3.2\Hz$ and an amplitude
ranging from $0$--$25\mm$. The peak response amplitude is limited to $\pm80\mm$;
at resonance, this limitation restricts the amplitude of the base motion to
approximately $\pm7\mm$. The motion of the base and the moving mass are measured
using laser displacement sensors (an Omron ZX2-LD100 and an Omron ZX1-LD300
respectively). In addition to the displacement measurements, the force provided
by the shaker is measured using an MCL-type load cell.

For the real-time control, a linear proportional-derivative (PD) controller is
used with manually tuned gains. The methodology is relatively insensitive to the
control gains used provided they are sufficient to stabilise any unstable orbits
that are encountered. The controller is implemented on a BeagleBone Black fitted
with a custom data acquisition board (hardware schematics and associated
software are open source and freely available~\cite{Barton2015rtc}). All
measurements are made at $1\kHz$ with no filtering.

A random perturbation signal is generated, when necessary, on the real-time
control board using the Box-Muller transformation to generate Gaussian
pseudo-random numbers which are then filtered using a sixth-order Butterworth
filter with a cut-off arbitrarily set at $10\Hz$ (below the Nyquist frequency of
the discretisation used for \eref{eqn:arx}).

Estimations of the Fourier coefficients of the response and the control action
are also calculated in real-time on the control board. However, this was for
convenience rather than a necessity.


\section{Experimental results}\label{sec:results}

\begin{figure}
  \centering
  \includegraphics{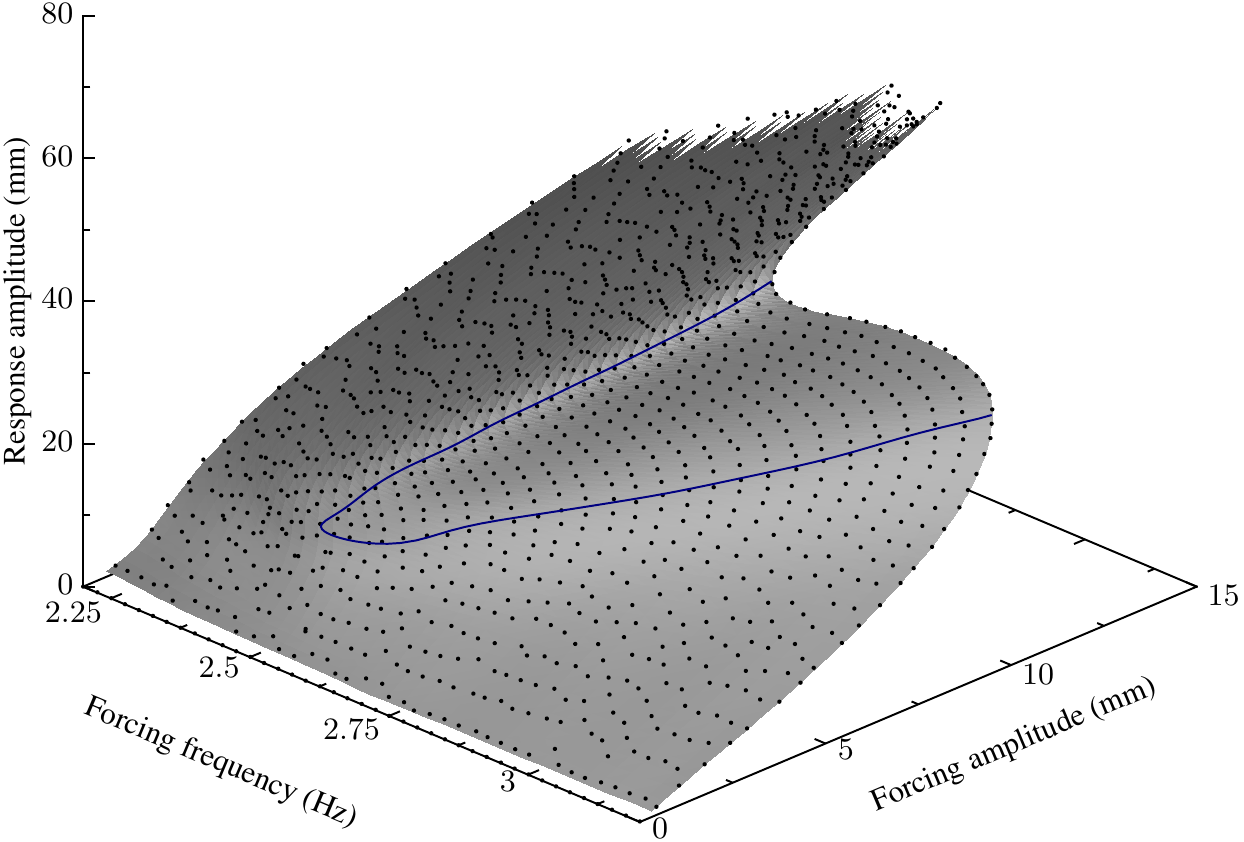}
  \caption{Measurements taken from $41$ continuation runs that vary the forcing
    amplitude for different fixed values of the forcing frequency. The forcing
    frequency is changed in steps between $2.2\Hz$ and $3.2\Hz$. The shaded
    surface is calculated using Gaussian process regression on the measured data
    points. The location of the unstable periodic orbits in this figure can be
    inferred from the geometry of the solution surface; the plotted curve (solid
    line) represents a 1D fold curve inside which are the unstable periodic
    orbits. This fold curve is calculated using the Gaussian process regressor.}
  \label{fig:solution_surface}
\end{figure}

The basic control-based continuation algorithm described in \sref{sec:cbc} was
used to do repeated continuations in the forcing amplitude (the amplitude of
displacement of the shaking table) for fixed values of the forcing
frequency. The forcing frequency, while fixed for individual continuation runs,
was varied between $2.2\Hz$ and $3.2\Hz$ in steps of $0.025\Hz$. At each data
point, full time series measurements were made. These are shown as black dots in
\fref{fig:solution_surface} where the forcing frequency and forcing amplitude
(in mm) are plotted against the response amplitude, which we define as the
magnitude of the first component in the Fourier series.

To aid visualisation, a continuous surface constructed from the individual data
points is also plotted in \fref{fig:solution_surface}. This continuous surface
is created using Gaussian progress regression on the collected data points where
the hyper-parameters for the Gaussian process are calculated by maximising the
marginal likelihood of the hyper-parameters~\cite[Sec. 5.4]{RasmussenWilliams2006}.

The use of Gaussian process regression (or any other similar scheme for
interpolating the multi-dimensional data) also allows for geometric features of
the solution surface to be easily extracted. One pertinent feature is the fold
in the solution surface which, from dynamical systems theory, indicates a change
in stability of the periodic orbits. As such, the fold curve shown in
\fref{fig:solution_surface} (blue curve) was extracted using numerical
continuation with \textsc{CoCo}~\cite{DankowiczSchilder2013} directly on the
regression surface defined by the Gaussian process.

\begin{figure}
  \centering
  \includegraphics{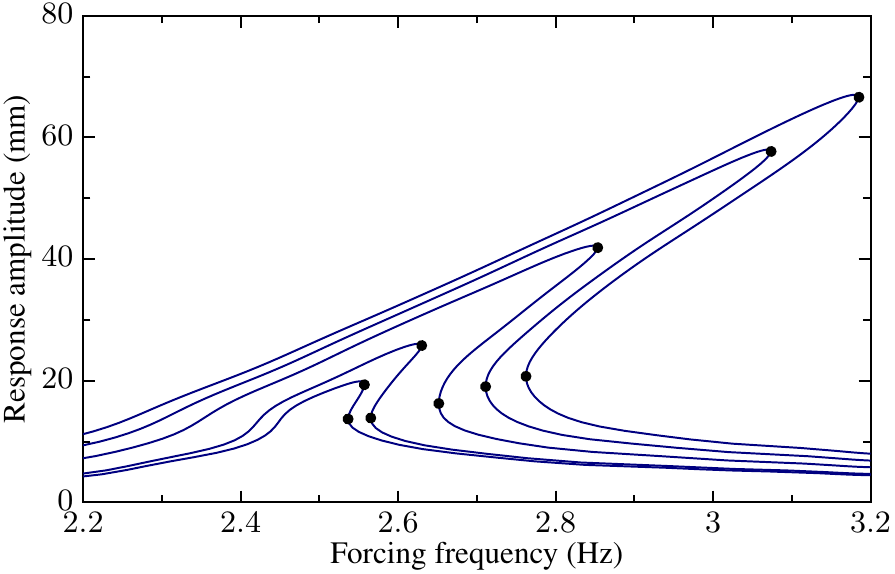}
  \caption{A series of frequency response curves extracted from the data shown
    in \fref{fig:solution_surface} using Gaussian process regression. The forcing
    amplitudes of the base motion are $1.9$, $2$, $2.5$, $3$ and $3.5\mm$
    respectively. Fold points (limit points) which determine the hysteresis
    region are marked with black dots.}
  \label{fig:freq_sweep}
\end{figure}

Other features of interest are frequency response curves which can also be
obtained through numerical continuation on the regression surface by fixing the
forcing amplitude to a prescribed value. \Fref{fig:freq_sweep} shows such
frequency response curves, including unstable periodic orbits, for fixed forcing
amplitudes of $\Gamma=1.9$, $2$, $2.5$, $3$ and $3.5\mm$. For high-amplitude
forcing, there is little to distinguish the results to those obtained from a
Duffing equation with hardening nonlinearity. However, for low-amplitude forcing
there seems to be a significant influence from frictional nonlinearities in the
bearing system suspending the mass.

\begin{figure}
  \centering
  \includegraphics{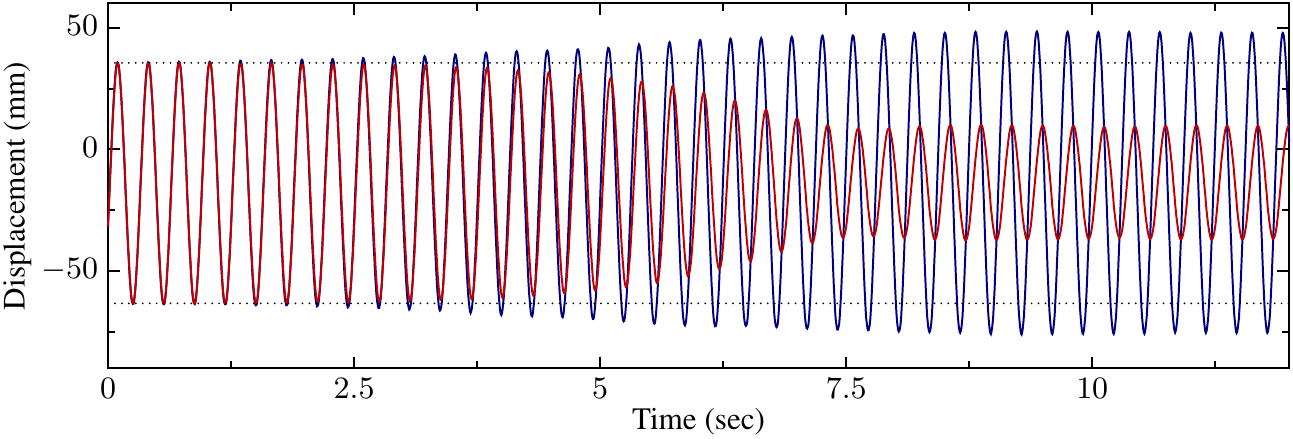}
  \caption{Two separate time-series measurements in taken open-loop conditions
    starting from the same unstable periodic orbit. The measurements where
    synchronised using the periodic forcing as a reference signal. Noise in the
    experiment randomly perturbs the trajectory from the unstable periodic orbit
    to either the stable high-amplitude orbit (blue) or the stable low-amplitude
    orbit (red). This shows that the stable manifold of the unstable periodic
    orbit acts as a separatrix between the basins of attraction of the high- and
    low-amplitude orbits.}
  \label{fig:unstable_escape_timeseries}
\end{figure}

In order to verify that the unstable orbits found in the experiment are true
unstable periodic orbits and not artefacts of the control scheme, we repeatedly
drive the system to a particular unstable periodic orbit and then turn off the
stabilising controller. As can be seen from the time series shown in
\fref{fig:unstable_escape_timeseries}, starting from the unstable periodic orbit
both the stable low-amplitude and stable high-amplitude periodic orbits can be
reached --- the stable manifold of the unstable orbit acts as a separatrix
between the two stable orbits as we would expect from dynamical systems
theory. Out of $40$ separate time series recorded starting from the same
unstable periodic orbit, $8$ end at the high-amplitude periodic orbit and the
remaining end at the low-amplitude periodic orbit.

\begin{figure}
  \centering
  \includegraphics{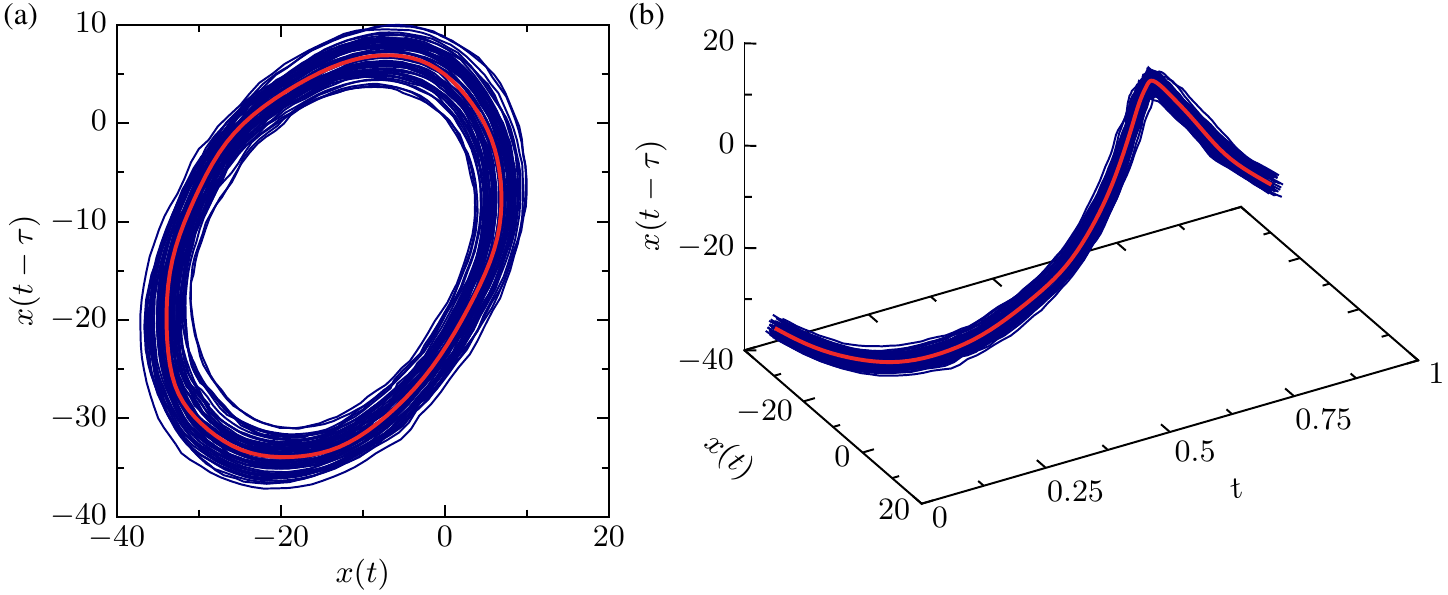}
  \caption{Panels (a) and (b) show 2D and 3D state-space projections of a single
    \newtext{stable} periodic orbit (red) and the randomly perturbed orbit
    (blue) used to calculate the stability of the periodic orbit. Time-delay
    coordinates are used as a proxy for the derivative of the position $x(t)$ to
    recreate the state space. The time coordinate in panel (b) is normalised
    such that a period of forcing takes one time unit.}
  \label{fig:orbit_delay_embedding}
\end{figure}

In order to apply the method outlined in \sref{sec:linear}, once an orbit has
been obtained using control-based continuation it must be perturbed with a
random input signal. One such example with a perturbation size of $0.5$ is shown
in \fref{fig:orbit_delay_embedding}. \newtext{(Strictly speaking, the
  perturbation size is measured in volts as it is an input to the shaker;
  however, the spectral content of the perturbation combined with the
  non-trivial frequency response of the shaker mean that it is not
  straightforward or useful to state the perturbation size in mm.)} To avoid
estimating the velocity of motion, a time-delay coordinate $x(t-\tau)$ is used
to reconstruct the state-space of the experiment. Here the value of $\tau$ used
is $T/5$ where $T$ is the period of the forcing.

\begin{figure}
  \centering
  \includegraphics{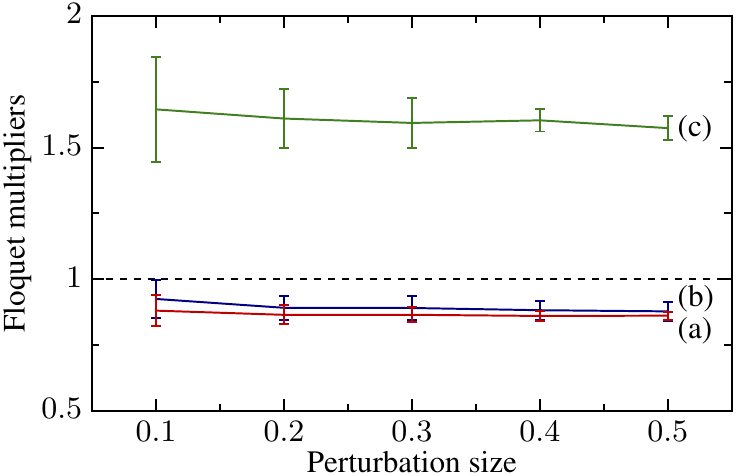}
  \caption{The mean and $95\%$ confidence interval of the absolute values of the
    \newtext{dominant} Floquet multipliers estimated against the size of the
    perturbation applied. Each point represents $10$ separate measurements. The
    points marked (a) (red) are calculated for a stable periodic orbit in
    open-loop conditions; the points marked (b) (blue) are calculated for the
    same stable periodic orbit in closed-loop conditions; and the points marked
    (c) (green) are calculated for an unstable periodic orbit in closed-loop
    conditions. \newtext{For the stable periodic orbits (a) and (b) the Floquet
      multipliers are complex; the estimated errors in the real and imaginary
      parts are approximately equal.}}
  \label{fig:control_effect}
\end{figure}

There is a single parameter in \sref{sec:linear} for which there is no
algorithmic way to determine an appropriate value, that is the amplitude of the
perturbation applied to the periodic orbit. Thus in order to determine an
appropriate amplitude we select two periodic orbits, one stable and one unstable
and calculate the Floquet multipliers using the fitted linear time-varying model
\eref{eqn:arx} for a variety of perturbation sizes and repeat the experiments
$10$ times. 

Throughout this paper we set $m=10$ and $n=4$ in \eref{eqn:arx}. 

The results of the Floquet multiplier estimations are shown in
\fref{fig:control_effect}. The points marked (a) (in red) are the absolute
values of the Floquet multipliers estimated for a stable periodic orbit while
the experiment is running in open-loop; similarly, the points marked (b) (in
blue) are the Floquet multipliers estimated for the same stable periodic orbit
while the experiment is running in closed-loop. Finally, the points marked (c)
(in green) are the Floquet multipliers estimated for an unstable periodic orbit
while the experiment is running in closed-loop. Each Floquet multiplier is
marked with the $95\%$ confidence range.

For both the stable and the unstable periodic orbits it can be seen that the
$95\%$ confidence range narrows significantly for larger perturbation sizes; the
larger perturbations allow the Floquet multipliers to be estimated more
consistently. However, the consistency of the results does not imply accuracy
--- as with estimating derivatives from finite differences, taking a large step
introduces errors caused by the nonlinearities in the system. Though, since the
mean magnitude of the Floquet multipliers does not change significantly beyond a
perturbation size of $0.2$, we can have reasonable confidence in the results and
so a perturbation size of $0.5$ is used throughout the remainder of this paper.

Furthermore, \fref{fig:control_effect} shows that running the system in
closed-loop rather than open-loop does not have a significant affect on the
estimation of the Floquet multipliers --- the error bars of the points (a) and
(b) overlap considerably.

Unfortunately no accurate independent estimations of the Floquet multipliers are
available to check the results of the system identification. Convergence tests
on the stable periodic orbit were performed, however the transient dynamics of
the electrodynamic shaker as the experiment equilibrated rendered the results
meaningless. Furthermore, escape tests from the unstable periodic orbit such as
those seen in \fref{fig:unstable_escape_timeseries} produced estimates ranging
from $1.2$ to $1.5$ depending on the measurement cut-offs used. As such, we
consider the best way to judge the accuracy of the estimations is via comparison
with geometric phenomena such as fold points where a Floquet multiplier should
pass through $+1$ in the complex plane.

\begin{figure}
  \centering
  \includegraphics{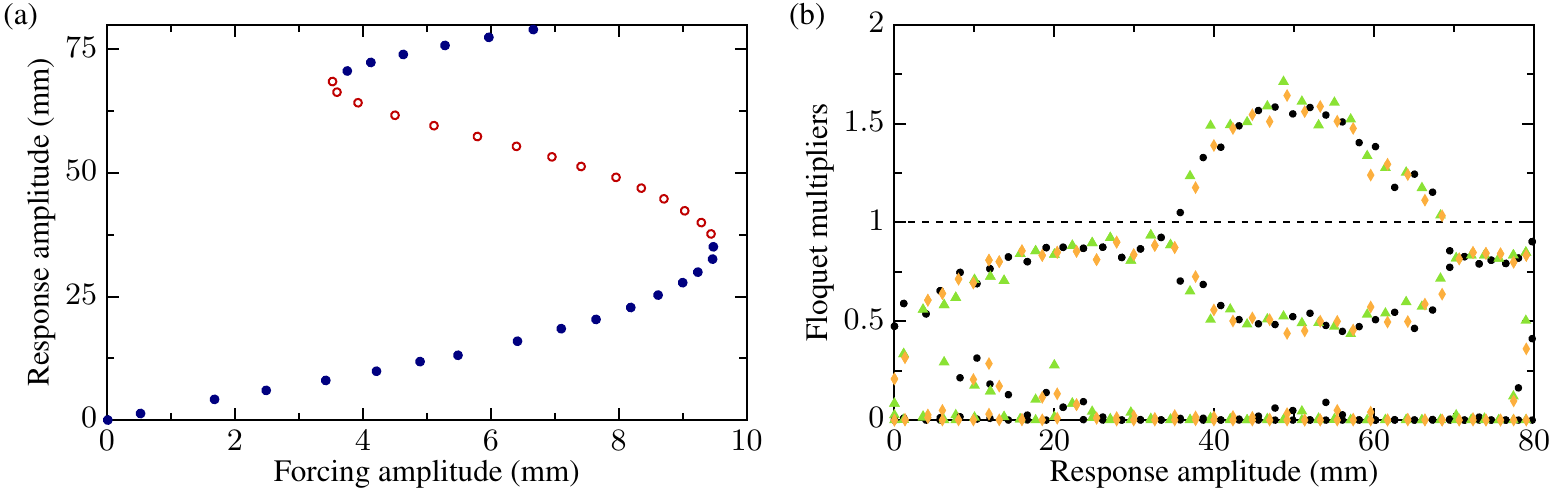}
  \caption{Panel (a) shows a single continuation \newtext{with the response
      amplitude plotted against the forcing amplitude} for a fixed fixed forcing
    frequency of $3.2\Hz$. Stable periodic orbits are shown by solid blue dots
    and unstable orbits are shown by red circles. The maximum displacement of
    the experiment is limited to $\pm80\mm$. Panel (b) shows the
    \newtext{absolute value of the estimated Floquet multipliers plotted against
      the response amplitude} for three separate continuation runs (each shown
    with a different symbol). It can be seen that most of the stable periodic
    orbits have complex conjugate Floquet multipliers which then become real
    close to the saddle-node bifurcations; between the two saddle-node
    bifurcations, one of these Floquet multipliers lies outside the unit circle
    while the other lies inside the unit circle.}
  \label{fig:branch_stability}
\end{figure}

Subsequently, a continuation in the forcing amplitude was performed with Floquet
multipliers estimated for each obtained solution; the results are shown in
\fref{fig:branch_stability}(a) with stable periodic orbits marked as solid blue
dots and unstable periodic orbits marked as red circles. The results agree very
well with what is expected from dynamical systems theory; the orbits between the
two fold points are all unstable.

To ensure repeatability, the same continuation run was carried out three times
and the absolute values of the corresponding Floquet multipliers are plotted in
\fref{fig:branch_stability}(b) with different symbols for the different runs. As
can be seen from \fref{fig:branch_stability}(b) the results are very consistent
with all the runs showing good quantitative as well as qualitative agreement. As
expected, a complex conjugate pair of Floquet multipliers becomes real close to
the fold point after which a single real multiplier passes through the unit
circle. This process reverses close to the second fold point.

\begin{figure}
  \centering
  \includegraphics{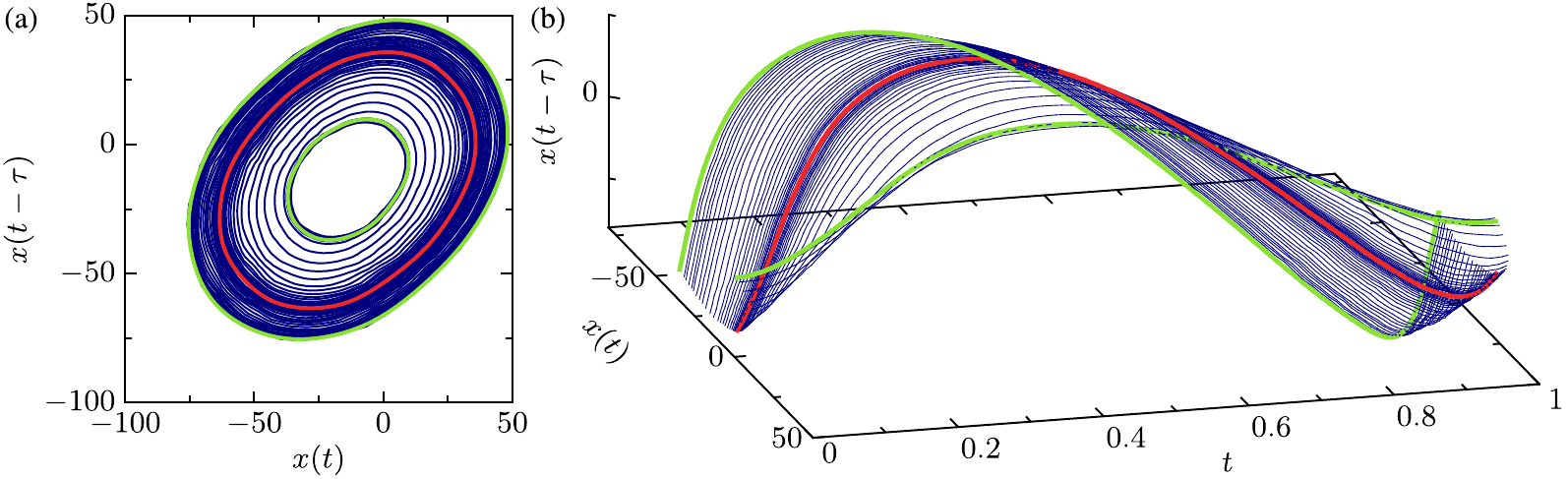}
  \caption{State-space projections in 2D and 3D of the data shown in
    \fref{fig:unstable_escape_timeseries} using time-delay coordinates. The red
    curve represents the unstable periodic orbit and the two green curves
    represent the stable high- and low-amplitude periodic orbits. The blue
    curves show the transient dynamics between the unstable orbits and the
    stable orbits. The two different stable orbits are reached with
    approximately equal probability. The time coordinate in panel (b) is
    normalised such that a period of forcing takes one time unit.}
  \label{fig:unstable_escape}
\end{figure}

As well as providing information about the Floquet multipliers of the periodic
orbits, the monodromy matrix obtained from \eref{eqn:arx} also provides the
stable and unstable eigendirections of the periodic
orbit. \Fref{fig:unstable_escape} shows the time series from
\fref{fig:unstable_escape_timeseries} plotted in state-space, with the unstable
periodic orbit (marked in red) lying between the two stable periodic orbits
(marked in green). From \fref{fig:unstable_escape}(b) it can be seen that the
escape from the unstable periodic orbit occurs along a two-dimensional unstable
manifold.

\begin{figure}
  \centering
  \includegraphics{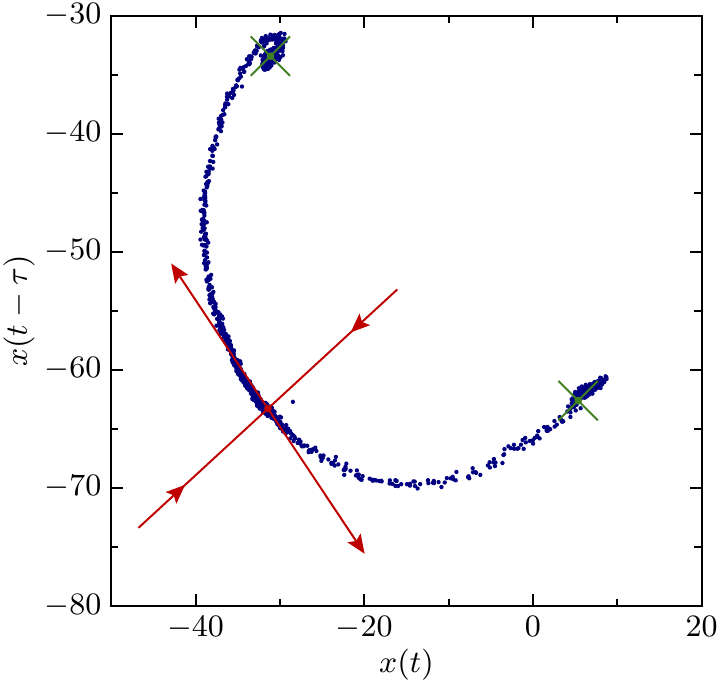}
  \caption{A Poincar\'e section showing measurements from $40$ separate time
    series of the escape from the unstable periodic orbit shown in
    \fref{fig:unstable_escape_timeseries} and \fref{fig:unstable_escape}. Each
    dot (blue) corresponds to an intersection with the Poincar\'e section
    defined by $t \mod T = 0$, where $T$ is the forcing period. The dots trace
    out the 1D unstable manifold of the unstable orbit. The positions of the
    stable high- and low-amplitude periodic orbits are marked with green
    crosses. The stable and unstable eigenspaces estimated using the methodology
    outlined in this paper are marked by red arrows. The unstable eigenspace
    shows good agreement with the dots from the time series measurements.}
  \label{fig:unstable_escape_poincare}
\end{figure}

\Fref{fig:unstable_escape_poincare} shows a Poincar\'e section at a fixed time
in the forcing cycle that was created with the data from
\fref{fig:unstable_escape} combined with data from a further $38$ independent
time series measurements; the blue dots denote intersections of the trajectories
with the Poincar\'e section. The intersection of the two-dimensional unstable
manifold with the Poincar\'e section results in a well defined one-dimensional
curve showing how trajectories leave the unstable periodic orbit (marked with a
red dot) and approach the two stable periodic orbits (marked with green crosses).

Superimposed on \fref{fig:unstable_escape_poincare} is a set of red arrows which
mark the unstable and stable eigendirections calculated from the monodromy
matrix of the periodic orbit using only local perturbations (as previously
described). The eigendirections represent the linearisation of the manifolds at
the unstable periodic orbit and, as such, show remarkably good agreement with
the measured unstable manifold which was calculated in open-loop
conditions. This opens up the possibility for using eigendirection information
in other calculations on the physical system, for example the estimation of
basins of attraction.


\section{Conclusions}\label{sec:conclusions}

This paper has proposed a new method for estimating Floquet multipliers and
their associated eigendirections for periodic orbits that encountered when using
control-based continuation on a physical experiment. A linear time-varying model
in the form of an auto-regressive model with exogenous inputs (ARX) is fitted to
each periodic orbit using small perturbations to the orbit to explore the nearby
state-space.

The method was demonstrated on a nonlinear mass-spring-damper-type experiment
that has a hardening spring characteristic. The nonlinearity is provided by
placing springs perpendicular to the direction of motion, thus creating a
geometric nonlinearity.  The Floquet multiplier estimations were shown to agree
with what is expected from dynamical systems theory and the associated
eigendirections match well with open-loop measurements taken.


\section{Data statement}\label{sec:data}

All the experimental data used in this paper has been deposited into the
University of Bristol Research Data Repository and is publicly available for
download~\cite{Barton2015stabdata}.


\section{Acknowledgements}\label{sec:ack}

D.A.W.B.\ is supported by EPSRC First Grant EP/K032739/1. D.A.W.B\ also
gratefully acknowledges discussions with and input from Ludovic Renson, Alicia
Gonzalez-Buelga and Simon Neild.


\section{References}\label{sec:refs}


\end{document}